\documentclass[leqno,10pt]{amsart}

\usepackage{amsmath}

\usepackage{amssymb}
\usepackage{graphicx}

\newtheorem{thm}{Theorem}

\newtheorem{lemma}[thm]{Lemma}

\theoremstyle{definition}
\newtheorem{defn}{Definition}

\theoremstyle{remark}
\newtheorem{remark}{Remark}

\numberwithin{equation}{section}

\def\R{\mathbb{R}}

\def\Z{\mathbb{Z}}
\def\E{\mathbb{E}}

\def\SO{\mathrm{SO}}

\def\I{\mathcal{I}}

\def\dim{\mathrm{dim\,}}

\begin{document}

\title[]{The geometric Cauchy problem for\\ the membrane shape equation}

\author{Gary R. Jensen}
\address{(G. R. Jensen)
Department of Mathematics, Washington University, One Brookings Drive,
St. Louis, MO 63130, USA}
\email{gary@math.wustl.edu}

\author{Emilio Musso}
\address{(E. Musso) Dipartimento di Scienze Matematiche, Politecnico di Torino,
Corso Duca degli Abruzzi 24, I-10129 Torino, Italy}
\email{emilio.musso@polito.it}

\author{Lorenzo Nicolodi}
\address{(L. Nicolodi) Di\-par\-ti\-men\-to di Ma\-te\-ma\-ti\-ca e Informatica,
Uni\-ver\-si\-t\`a degli Studi di Parma, Parco Area delle Scienze 53/A,
I-43124 Parma, Italy}
\email{lorenzo.nicolodi@unipr.it}

\thanks{Authors partially supported by
PRIN 2010-2011 ``Variet\`a reali e complesse: geometria, to\-po\-lo\-gia e analisi ar\-mo\-ni\-ca'';
FIRB 2008 ``Geometria Differenziale Complessa e Dinamica Olomorfa'';
and by the GNSAGA of INDAM}


%\date{Version of September 30, 2014}
%\date{Version of \today}

%\dedicatory{}


\begin{abstract}
We address the geometric Cauchy problem for surfaces associated to the {membrane shape equation}
describing equilibrium configurations of vesicles formed by lipid bilayers.
This is the Euler--Lagrange equation of the Canham--Helfrich--Evans elastic curvature energy
subject to constraints on the enclosed volume and the surface area.
Our approach uses the method of moving frames and techniques from the theory of
exterior differential systems.
\end{abstract}

\maketitle

\section{Introduction}\label{s:intro}

Lipid bilayers are the basic elements of biological membranes and constitute
the main separating structure of living cells. In aqueous solution, lipid bilayers typically form
closed surfaces or vesicles (closed bilayer films) which exhibit a large variety of different shapes,
such as the non-spherical, biconcave shape characteristic of red blood cells
\cite{DH1976, Lip, Seif}.
Since in most biologically significant cases the width of the membrane exceeds the thickness
by several orders of magnitude, the membrane can be regarded as
a two-dimensional surface $S$ embedded in three-dimensional space $\R^3$, with
enclosed volume $V(S)$ and surface area $A(S)$.

\vskip0.1cm
A continuum-mechanical description of lipid bilayers
which has been shown to be predictive with respect to the observed
shapes dates from the work of Canham \cite{Can}, Helfrich \cite{Hel}, and Evans \cite{Ev}
at the beginning of the 1970s.
In the Canham--Helfrich--Evans model, an equilibrium configuration $S$,
at fixed volume and surface area, is determined by minimization of the elastic bending
energy \cite{Hel, HKS1992, Ni1993, OY-L-X-book1999, Pet1985, Steig1999, Steig2003, TOY2004, vHemm}
\begin{equation}\label{helfric}
 \mathcal F(S) =
 \frac{k}{2}\int_S{(2H + c_0)^2dA} + \bar k\int_S{KdA}
 + p {V(S)} + \lambda {A(S)}.
  \end{equation}
%
% in \cite{OYH1989}, the first term of the functional is
% $\frac{k}{2}\int_S{(c_1+c_2 - c_0)^2dA}$; note that in \cite{OYH1989}, p. 5281,
% $H = (1/2) g^{ij}L_{ij} = (-c_1 -c_2)/2$; note that, while the usual def of $H$ (Stoker,
% Hopf, do Carmo) is used, a different sign for the principal curvatures is taken
%
Here $H = (a +c)/2$ and $K = ac$ are the mean and Gauss curvatures of the surface $S$,
where $a$ and $c$ denote the principal curvatures,
$dA$ is the area element,
and $k$,
$\bar k$,
 $c_0$, $p$, and $\lambda$ are constants:
$k$, $\bar k$ are material dependent elasticity parameters,
$c_0$ is the spontaneous curvature,
$p$ denotes the difference between the outside and the inside pressure,
and $\lambda$ is the surface lateral tension. The pressure $p$ and the tension $\lambda$
play the role of Lagrange multipliers for the constraints of constant volume and
surface area.

\vskip0.1cm
The Euler--Lagrange equation of the functional \eqref{helfric}, %which is
often referred to as the
{\em membrane shape equation} or the {\em Ou-Yang--Helfrich equation},
takes the form \cite{OYH1987, OYH1989}
\begin{equation}\label{shape-eq-alt1}
 2k\left\{\Delta H + 2H(H^2 -K)\right\} - \left(2\lambda + k c_0^2\right)H
  - 2k c_0 K + p = 0,
  \end{equation}
where $\Delta$ denotes the Laplace--Beltrami operator of the induced metric on $S$.
This implies that an equilibrium shape surface $S$
satisfies a fourth order nonlinear partial differential equation of the form
\begin{equation}\label{shape-eq}
\Delta H  = \Phi(a,c),
\end{equation}
where $\Phi$ is a real analytic symmetric function of the principal curvatures $a$,
$c$. Examples of surfaces associated to \eqref{shape-eq-alt1} include
Willmore surfaces \cite{Blaschke, Br-duality, Thomsen, Willmore},
which are solutions of the differential
equation
\begin{equation}\label{willmore-eq}
\Delta H + 2H(H^2 -K) = 0.
\end{equation}
Willmore surfaces are best known in connection to the celebrated Willmore conjecture,
recently confirmed by Marques and Neves \cite{Ma-Ne}.

\vskip0.1cm
The aim of the present paper is to solve the geometric
Cauchy problem for the class of surfaces in $\R^3$ satisfying
the membrane shape equation. More precisely, we shall prove the following.

\begin{thm}\label{main-thm}
Let $\alpha : J
%=(\mathrm a,\mathrm b)
\to \R^3$ be a
real analytic curve with $\|\alpha'(x) \|=1$,
where $J\subset \R$ is an open interval.
Let $(T = \alpha'$, $N$, $B)$ be the Frenet frame field
along $\alpha$, with Frenet--Serret equations
\[
 T' = \kappa N, \quad N' = -\kappa T + \tau B, \quad B' = -\tau N,
\]
where $\kappa(x)\neq 0$ is the curvature and $\tau(x)$ is the torsion
of $\alpha$.
For $x_0\in J$, consider the
unit normal vector
$W_0 =  N(x_0) \cos{a_0} + B(x_0)\sin{a_0}$.
Let $h$,
$h^W : J \to \R$ be two real analytic functions and assume that
\begin{equation}\label{kt-ineq}
  h + \kappa\sin {\left(-\int_{x_0}^x\tau (u)du +a_0\right)} < 0.
   \end{equation}
Then, there exists a real analytic immersion $f : \Sigma \to \R^3$, where $\Sigma \subset \R^2$
is an open neighborhood of $J\times \{0\}$, with principal curvature line coordinates $(x,y)$,
such that:
\begin{enumerate}
\item the mean curvature $H$ of $f$ satisfies
\[
 \Delta H  = \Phi(a,c);
  \]

\item the restriction $f\vert_J = \alpha$;

\item the tangent plane to $f$ at $f(x_0, 0)$ is spanned by
$T(x_0)$ and $W_0$;

\item $\alpha$ is a curvature line of $f$;

\item  $H\vert_J = h$ and ${\frac{\partial H}{\partial y}}\big\vert_J = h^W$.

\end{enumerate}

\noindent Moreover, if $\hat f : \hat \Sigma \to \R^3$ is any other principal immersion
satisfying the above conditions, then
$f(\Sigma \cap \hat \Sigma) = \hat f (\Sigma\cap \hat \Sigma)$.

\end{thm}

The Cauchy problem addressed in Theorem \ref{main-thm}
can be viewed as a generalization of similar problems for constant mean curvature surfaces
in $\R^3$ and for constant mean curvature one surfaces in hyperbolic 3-space \cite{BD, GM},
which in turn are both inspired by the classical Bj\"orling problem for minimal surfaces in
$\R^3$ \cite{Ni1989}. Recently, the geometric Cauchy problem has been investigated for
several surface classes and in different geometric situations (see, for instance,
\cite{ACM2003, ACG2007, BS2013} and the references therein).
\vskip0.1cm
However, our approach to the problem is different
in that we use techniques from the Cartan--K\"ahler
theory of Pfaffian differential systems and the method of moving frames
(see \cite{MNjmp} and \cite{MNcag, MNphysD, MNima} for a similar approach to the integrable
system of Lie-minimal surfaces and other systems in submanifold geometry).
This work was mainly motivated by a paper of Tu and Ou-Yang \cite{TOY2004},
in which the authors propose a geometric scheme to discuss
the questions of the shape and stabilities of biomembranes
within the framework of exterior differential forms.
\vskip0.1cm
The first step in our discussion
consists in the construction of a Pfaffian differential system (PDS) whose integral manifolds
are canonical lifts of principal frames along surfaces satisfying the membrane shape equation \eqref{shape-eq}. We then compute the algebraic generators of degree two for such a PDS and show that the polar space of a 1-dimensional integral element is 2-dimensional. Next, by the Cartan--K\"ahler
theorem we deduce the existence of a unique real analytic integral surface
passing through a real analytic integral curve. Finally, we build 1-dimensional
integral curves from arbitrary real analytic space curves and two real analytic
functions.
The proof of Theorem \ref{main-thm} follows from a suitable geometric interpretation
of these results. Interestingly enough, in the proof of Theorem \ref{main-thm}, a crucial role
is played by the use of a
{\em relatively parallel
adapted frame} of Bishop \cite{Bishop} along the given curve $\alpha$.

\vskip0.1cm
The paper is organized as follows.
Section \ref{s:background} introduces some background material and recalls the basic facts
about Pfaffian differential systems in two independent variables.
Section  \ref{s:PDS} constructs the PDS
for the class of surfaces associated to the membrane shape equation.
Section \ref{s:PDS-involution} proves that this PDS is in involution.
Section \ref{s:main-thm} proves Theorem \ref{main-thm}.
Finally, Section \ref{s:ex} discusses some examples.

\vskip0.1cm
For the subject of exterior differential systems we refer the
reader to the monographs \cite{BCGGG, Cartan-book, Grbook, GJbook, ILlibro}.
The summation convention over repeated indices is used throughout the paper.

\vskip0.1cm
The authors would like to thank the referees for their useful comments and suggestions.

\section{Background material}\label{s:background}

\subsection{The Euclidean group and the structure equations}

Let $\E(3) = \R^3\rtimes \SO(3)$ be the Euclidean group of proper
rigid motions of $\R^3$.
A group element of $\E(3)$ is an ordered pair $(P,A)$, where $P\in \R^3$
and $A$ is a $3\times 3$ orthogonal matrix with determinant one. If we let
$A_j\in \R^3$, $j = 1,2,3$, denote the $j$-th column vector of $A$
and regard $P$ and the $A_j$ as $\R^3$ valued functions on $\E(3)$, there
exist unique left invariant 1-forms $\theta^i$ and $\theta^i_j$, $i,j = 1,2,3$,
such that
\begin{equation}\label{dP-dA}
dP = \theta^i A_i, \quad dA_j = \theta^i_j A_i, \quad j = 1,2,3.
\end{equation}
The 1-forms $\theta^i$, $\theta^i_j$ are the Maurer--Cartan forms of $\E(3)$.
Differentiating the orthogonality condition $A_i \cdot A_j = \delta_{ij}$ yields
$\theta^i_j = -\theta^j_i$, $i,j = 1,2,3$. These are the only relations among
the Maurer--Cartan forms and then $(\theta^1, \theta^2, \theta^3, \theta^2_1, \theta^3_1,
\theta^3_2)$ is a basis for the space of left invariant 1-forms on $\E(3)$.
Differentiating \eqref{dP-dA}, we obtain the structure equations of $\E(3)$
\begin{equation}\label{dthetai}
\begin{cases}
%\left\{\!\!\begin{array}{l}
 d\theta^1 = \theta^2_1 \wedge \theta^2 + \theta^3_1 \wedge\theta^3,\\
  d\theta^2 = -\theta^2_1 \wedge \theta^1 + \theta^3_2 \wedge\theta^3,\\
   d\theta^3 = -\theta^3_1 \wedge \theta^1 - \theta^3_2 \wedge\theta^2,
%    \end{array}\right.
\end{cases}
\end{equation}
and
\begin{equation}\label{dthetaij}
\begin{cases}
%\left\{\!\!\begin{array}{l}
 d\theta^2_1 = \theta^3_2 \wedge \theta^3_1,\\
  d\theta^3_1 = -\theta^3_2 \wedge \theta^2_1,\\
   d\theta^3_2 = \theta^3_1 \wedge \theta^2_1.
%    \end{array}\right.
\end{cases}
\end{equation}

\subsection{Principal frames and invariants}

Let $f : X \to \R^3$ be a smooth immersion of a connected, orientable
2-dimensional manifold $X$, with unit normal vector field $n$.
Consider the orientation of $X$ induced by $n$ from the orientation
of $\R^3$. Suppose that $f$ is free of umbilic points and
with the given normal vector such that the principal curvatures
$a$ and $c$ satisfy $a > c$.
A principal frame field along $f$ is a map $(f,A) : U \to \E(3)$
defined on some open connected set $U\subset X$, such that, for each $\zeta\in U$,
the tangent space
$df(T_\zeta X )=\text{span}\{{A_1}(\zeta), A_2(\zeta)\}$, $A_3(\zeta) = n(\zeta)$,
and ${A_1}(\zeta)$, $A_2(\zeta)$
are along the principal directions corresponding to $a(\zeta)$ and $c(\zeta)$,
respectively. Any other principal frame field on $U$ is of the form
$(f, (\pm A_1, \pm A_2, A_3))$. Thus, if $X$ is simply connected, possibly passing
to a double cover, we may assume the existence of a globally defined principal
frame field $(f,A)$ along $f$. Following the usual practice in the method of
moving frames we use the same notation to denote the forms on $\E(3)$ and
their pullbacks via $(f,A)$ on $X$.

Let $(f,A)$ be a globally defined principal frame field along $f$. Then on $X$,
$(\theta^1, \theta^2)$ defines a coframe field, $\theta^3$
vanishes identically, and
\[
  \theta^3_1 = a \theta^1,\quad \theta^3_2 = c \theta^2,
    \quad \theta^2_1 = p\theta^1 + q \theta^2,
       \]
where $a>c$ are the principal curvatures and $p, q$ are smooth functions,
the Christoffel symbols of $f$ with respect to $(\theta^1, \theta^2)$.
The structure equations of $\E(3)$ give
\begin{equation}\label{dtheta1-2}
 d\theta^1 = p \theta^1\wedge \theta^2,
\quad d\theta^2 = q \theta^1 \wedge \theta^2,
\end{equation}
the Gauss equation
\begin{equation}\label{gauss-eq}
 dp \wedge \theta^1 +dq \wedge \theta^2
+ (ac +p^2 +q^2) \theta^1\wedge \theta^2 =0,
\end{equation}
and the Codazzi equations
\begin{equation}\label{codazzi-eq}
\begin{cases}
 da \wedge \theta^1  +  p(c-a) \theta^2 \wedge \theta^1 =0,\\
 dc \wedge \theta^2  +  q(c-a) \theta^1 \wedge \theta^2 =0.
\end{cases}
\end{equation}

With respect to the principal frame field $(f,A)$, the first and second fundamental
forms of $f$ are given by
\[
  I = df \cdot df = \theta^1\theta^1 + \theta^2\theta^2,
\quad II = - df \cdot dn = a\theta^1\theta^1 +c \theta^2\theta^2.
\]
For any smooth function $g : X \to \R$, we set
\begin{equation}\label{i}
  dg = g_1\theta^1 +  g_2 \theta^2,
   \end{equation}
where the functions $g_1, g_2 : X \to \R$ play the role of partial derivatives relative
to $\theta^1, \theta^2$. In general, mixed partials are not equal, but satisfy
\begin{equation}\label{ii}
   g_{12} - g_{21} = p g_1 + q g_2,
   \end{equation}
by \eqref{dtheta1-2}.
Using the relation
\[
 (\Delta g )\,\theta^1\wedge \theta^2=d\ast d g
   =d(- g_2\theta^1 + g_1\theta^2),
    \]
which defines
the Laplace--Beltrami operator $\Delta$ of the metric $I$ on $X$, we find
\[
 \Delta g =
  g_{11} + g_{22} + q g_1 - p g_2.
   \]

In this formalism, the Gauss and Codazzi equations can be written as
\begin{equation}\label{gauss-codazzi-eq}
\begin{cases}
 p_2 - q_1 = ac +p^2 +q^2, \\
 a_2 = - p(c-a),\\
  c_1 =  - q(c-a).
\end{cases}
\end{equation}
From this we see that there exists a smooth function $r$, such that
\begin{equation}\label{r}
\begin{cases}
  p_2 = r + \frac{1}{2}(ac +p^2 +q^2),\\
   q_1 = r - \frac{1}{2}(ac +p^2 +q^2).
    \end{cases}
\end{equation}
Differentiating the structure equations \eqref{gauss-codazzi-eq} and using \eqref{i} and
\eqref{ii}, we get
\begin{equation}\label{gauss-codazzi-eq-r}
\begin{cases}
a_{21} = (p_1-pq) (a-c) + p a_1, \quad a_{12} = 2pa_1 +p_1(a-c), \\
 a_{22}= \left( r + \frac{1}{2}(ac +p^2 +q^2) \right)(a-c) + p^2(a-c) -pc_2,\\
  c_{21} = q_2 (a-c)  - 2q c_2 , \quad c_{12} = (q_2 + pq)(a - c) - q c_2, \\
   c_{11}= \left( r - \frac{1}{2}(ac +p^2 +q^2) \right)(a-c) + qa_1 - q^2(a-c).\\
\end{cases}
\end{equation}

Applying these formulae, we are now able to express the Laplace--Beltrami operator of
the mean curvature of the immersion $f$ in terms of the {functions}
$p, q, a, c, r, a_1, c_2$,
\begin{equation}\label{delta-H}
\Delta H =
 \frac{1}{2}\left( a_{11} +  c_{22}\right)
  - r (c-a) + q a_1 -p c_2.
\end{equation}
Thus, $f$ satisfies the differential relation
$\Delta H = \Phi(a,c)$ if and only if
\begin{equation}\label{psi}
 a_{11} +  c_{22}
%=\Psi (p,q,a,c,r, a_1, c_2),
= 2\left(\Psi(p,q,a,c,a_1,c_2) + r(c-a) \right)
\end{equation}
where
\[
\Psi (p,q,a,c, a_1, c_2)
%=2\left[\Phi(a,c) + r (c-a) + p c_2 - q a_1\right].
= \Phi(a,c) + pc_2 - qa_1.
\]

%Remark: the point is that \Psi should not be a function of r.

\subsection{Pfaffian differential systems in two independent variables}
In this section we recall some basic facts about the Cartan--K\"ahler
theory of Pfaffian differential systems in two independent variables.
Let $M$ be a smooth manifold
and let
\[
 \left(\omega^1, \omega^2, \eta^1, \dots, \eta^k, \pi^1, \dots, \pi^s \right)
  \]
be a coframe field on $M$. Consider the 2-form $\Omega = \omega^1 \wedge \omega^2$,
which is referred to as the {\em independence condition}, and let
\[
  \I = \{\eta^1, \dots, \eta^k, d\eta^1, \dots, d\eta^k   \}
  \]
be the ideal of the algebra
of exterior differential forms on $M$ generated by the 1-forms
$\eta^1, \dots, \eta^k$ and the 2-forms $d\eta^1, \dots, d\eta^k$.
The pair $(\I, \Omega)$ is called a {\em Pfaffian differential system
in two independent variables}. A 2-dimensional integral manifold (integral surface) of $(\I,\Omega)$
is an immersed 2-dimensional manifold $X\subset M$, not necessarily embedded, such that
\[
  \Omega\vert_X \neq 0, \quad \eta^j\vert_X = d\eta^j\vert_X =0, \quad j=1, \dots, k.
  \]
Similarly, a 1-dimensional integral manifold (integral curve) of the Pfaffian system
$(\I, \Omega)$ is an immersed 1-dimensional manifold
%curve
$\mathcal A \subset M$, such that
\[
  \omega^1\omega^1 + \omega^2\omega^2\vert_{\mathcal A} \neq 0, \quad
   \eta^j\vert_{\mathcal A} =0, \quad j=1, \dots,k.
    \]
A 1-dimensional integral element
of $(\I,\Omega)$, at a fixed point $m\in M$, is a 1-dimensional subspace $[\xi]$ of $T_mM$,
spanned by a non-zero tangent vector $\xi$, such that
\[
  \omega^1(\xi)\omega^1(\xi) + \omega^2(\xi)\omega^2(\xi) \neq 0, \quad
   \eta^j(\xi) = 0, \quad j= 1,\dots, k.
    \]
Similarly, a 2-dimensional integral element of $(\I,\Omega)$
is a 2-dimensional subspace $W$ of $T_mM$, such that
\[
  \Omega\vert_W \neq 0, \quad \eta^j\vert_W = d\eta^j\vert_W =0, \quad j=1, \dots, k.
  \]
Given a 1-dimensional integral element $[\xi]$, its {\em polar space} $H([\xi])$ is the
subspace of $T_mM$ defined by the linear equations:
\[
   \eta^j = 0, \quad  i_\xi   d\eta^j = 0, \quad j = 1, \dots, k.
   \]
If
\begin{equation}\label{dim-polar-space}
  \dim H([\xi]) = 2, \quad \text{and} \quad \Omega\vert_{H([\xi])} \neq 0,
   \end{equation}
then $H([\xi])$ is the unique 2-dimensional integral element containing $[\xi]$.
The following statement is a consequence of the general Cartan--K\"ahler theorem
for exterior differential systems in involution \cite{BCGGG, Cartan-book}.

\begin{thm}
Let $M$ be a real analytic manifold and assume the the coframe field
\[
 \left(\omega^1, \omega^2, \eta^1, \dots, \eta^k, \pi^1, \dots, \pi^s \right)
  \]
is real analytic. If condition \eqref{dim-polar-space} is satisfied for every
1-dimensional integral element, then for any real analytic 1-dimensional integral
curve $\mathcal A \subset M$ there exists a 2-dimensional integral manifold
$X\subset M$ such that $\mathcal A\subset X$. The manifold $X$ is unique, in the
sense that if $\tilde X$ is another 2-dimensional integral manifold passing
through $\mathcal A$, then there exists an open neighborhood $U\subset M$
of $\mathcal A$ such that $\tilde X\cap U = X \cap U$.
\end{thm}

\section{The Pfaffian differential system of lipid bilayer membranes}\label{s:PDS}

\subsection{Principal frames}

Let $Y_{(1)}$ be the 10-dimensional manifold
\[
  Y_{(1)} = \E(3) \times \left\{(p,q,a,c)\in \R^4 \,| \, a -c > 0 \right\}.
    \]
Let
\[
  \theta^1, \theta^2, \theta^3, \theta^2_1, \theta^3_1, \theta^3_2
  \]
be the basis of left invariant Maurer--Cartan forms of $\E(3)$ which satisfy structure equations
\eqref{dthetai} and \eqref{dthetaij}.

On $Y_{(1)}$ consider the Pfaffian differential system $(\I_1, \Omega)$
differentially generated by the 1-forms
\begin{equation}\label{alphas}
 \begin{array}{ll}
  \alpha^1 = \theta^3, &
   \alpha^2 = \theta^2_1 -p\theta^1-q\theta^2,\\
    \alpha^3 = \theta^3_1 -a\theta^1, &
     \alpha^4 = \theta^3_2 -c\theta^2,\\
       \end{array}
        \end{equation}
with independence condition $\Omega = \theta^1 \wedge\theta^2$.

\begin{remark}
The integral manifolds of $(\I_1, \Omega)$ are the {\it second order prolongations of umbilic free
immersed surfaces} in $\R^3$, i.e., smooth maps
\[
   F_{(1)} := (F,A,p,q,a,c) : X \to Y_{(1)}
   \]
defined on an oriented, connected 2-dimensional manifold $X$ such that:

\begin{itemize}
\item  $F : X \to \R^3$ is an umbilic free smooth immersion;

\item $(F,A)= (F, (A_1,A_2,A_3)) : X \to \E(3)$ is a principal frame field along $F$;

\item $(\theta^1,\theta^2)$ is a positively oriented orthonormal coframe, dual to the
trivialization of $dF(T(X))$ defined by the tangent vector fields $A_1$, $A_2$
along $F$;

\item $\theta^1\theta^1 + \theta^2\theta^2$ is the first fundamental form
induced by $F$;

\item $A_1, A_2$ are tangent to the principal curvature
lines of $F$;

\item $\theta^2_1 = p\theta^1 +q\theta^2$ is the Levi Civita connection with
respect to $(\theta^1, \theta^2)$;

\item $a$ and $c$ are the principal curvatures of $F$ and $a >c$;

\item $a\theta^1\theta^1 + c \theta^2\theta^2$ is the second fundamental
form of $F$;

\item $A_3$ is the Gauss map of $F$.

\end{itemize}
Note that $F_{(1)}$ is uniquely determined by the orientation of $X$, by the
umbilic free immersion $F$, and by the assumption that $a > c$.

\end{remark}

From the structure equations \eqref{dthetai} and \eqref{dthetaij}, we obtain,
modulo the algebraic ideal $\{\alpha^j \}$ generated by the 1-forms
$\alpha^1, \dots, \alpha^4$,
\begin{equation}\label{eq-Y1-1}
\begin{cases}
d\theta^1 \equiv p \theta^1 \wedge \theta^2,\\
 d\theta^2 \equiv q \theta^1 \wedge \theta^2,\\
\end{cases}
\mod \{\alpha^j \}
\end{equation}
and
\begin{equation}\label{eq-Y1-2}
\begin{cases}
 d\alpha^1 \equiv 0,\\
  d\alpha^2 \equiv - dp\wedge\theta^1 -dq \wedge \theta^2
      - (ac +p^2 +q^2)\theta^1\wedge\theta^2,\\
   d\alpha^3 \equiv - da\wedge\theta^1 +p(c-a)\theta^1 \wedge \theta^2,\\
    d\alpha^4 \equiv - dc\wedge\theta^2 - q(c-a)\theta^1 \wedge \theta^2.\\
\end{cases}
\mod \{\alpha^j \}
\end{equation}

\subsection{Prolongation}

Let $Y_{(2)}$ be the 17-dimensional manifold
\[
  Y_{(2)} = Y_{(1)}\times \left\{(p_1, q_2, r,  a_1, c_2,  a_{11}, c_{22})\in \R^7\right\},
   \]
where
\[
  (p_1, q_2, r,  a_1, c_2,  a_{11}, c_{22})
  \]
are the new fiber coordinates. Consider on $Y_{(2)}$ the Pfaffian differential
system $(\I_2, \Omega)$ differentially generated by the 1-forms
\[
  ( \alpha^1, \dots, \alpha^4, \beta^1, \beta^2,
   \gamma^1, \gamma^2, \delta^1, \delta^2),
    \]
where $\alpha^1, \dots, \alpha^4$ are defined as in \eqref{alphas} and
\begin{equation}\label{betas-gammas-deltas}
\begin{cases}
\beta^1 = dp - p_1 \theta^1 - (r+\frac{1}{2}(ac +p^2 +q^2))\theta^2,\\
  \beta^2 = dq - (r- \frac{1}{2}(ac +p^2 +q^2))\theta^1 - q_2 \theta^2,\\
   \gamma^1 = da - a_1 \theta^1 + p(c-a) \theta^2, \\
    \gamma^2 = dc + q(c-a) \theta^1  - c_2 \theta^2, \\
     \delta^1 =  d a_1 - a_{11}\theta^1 +
       (-2 a_1 p +(c-a)  p_1)\theta^2, \\
      \delta^2 =  d c_2 + (2 c_2 q +(c-a) q_2)\theta^1
        -  c_{22}\theta^2. \\
\end{cases}
\end{equation}
From \eqref{eq-Y1-2} and \eqref{betas-gammas-deltas} it follows that,
modulo the algebraic ideal generated by $\alpha^j$, $\beta^a$, $\gamma^a$ and
$\delta^a$, $j=1,2,3,4; a = 1,2$, we have
\begin{equation}\label{d-alphas-d-gammas=0}
\begin{cases}
  d\alpha^j \equiv 0, \quad j=1,2,3,4, \\
    d\gamma^a \equiv 0, \quad a=1,2, \\
\end{cases}
\mod \{\alpha^j, \beta^a, \gamma^a, \delta^a \}
\end{equation}

\begin{remark}
The integral manifolds of $(\I_2, \Omega)$ are {\it fourth order prolongations
of immersed surfaces} in $\R^3$, that is, smooth maps
\[
  F_{(2)} := (F, A, p,q,a,c, p_1, q_2, r,  a_1, c_2,  a_{11}, c_{22}) : X \to Y_{(2)},
     \]
where $(F,A,p,q,a,c) : X \to Y_{(1)}$ is an extended principal frame along $F$
and the functions
\[
  p_1, q_2, r,  a_1, c_2,  a_{11}, c_{22} : X \to \R
     \]
are defined by
\[
\begin{cases}
  dp = p_1\theta^1 + p_2\theta^2, \quad dq = q_1\theta^1 + q_2\theta^2, \\
     r = p_2 - \frac{1}{2}(ac + p^2 +q^2)= q_1 + \frac{1}{2}(ac + p^2 +q^2), \\
     da = a_{1}\theta^1 + a_2\theta^2, \quad dc = c_{1}\theta^1 + c_2 \theta^2,\\
     da_1 = a_{11}\theta^1 + a_{12}\theta^2,\quad dc_2 = c_{21}\theta^1 + c_{22}\theta^2.
     \end{cases}
     \]
Note that $F_{(2)}$ is uniquely determined by the orientation of $X$, by
the immersion $F$, and by the assumption that the principal curvature
satisfy $a > c$.

The canonical prolongations of smooth immersions
$F : X \to \R^3$ satisfying the partial differential relation
\[
 \Delta H = \Phi(a,c)
  \]
are characterized by the equation
\begin{equation}
 a_{11} +   c_{22} = 2\left[\Phi(a,c) + r (c-a)
   -q a_1 + p  c_2 \right].
    \end{equation}
\end{remark}

\subsection{The Pfaffian differential system of lipid bilayer membranes}

Let $Y_{\ast}$ be the 16-dimensional real analytic submanifold of $Y_{(2)}$ defined by
\[
   a_{11} +  c_{22} = 2\left[\Phi(a,c) + r (c-a)
   -q a_1 + p   c_2 \right].
     \]
On $Y_{\ast}$ we consider the fiber coordinates
\[
 p, q, a, c,  p_1,   q_2,  r,   a_1,   c_2,
       \ell,
    \]
where $\ell$ is given by
\begin{equation}\label{ell}
\begin{cases}
 a_{11} = \ell +r (c -a) + \Psi(p,q,a,c, a_1, c_2), \\
  c_{22} = -\ell +r (c -a) + \Psi(p,q,a,c, a_1, c_2), \\
   \end{cases}
     \end{equation}
where
\[
  \Psi(p,q,a,c, a_1, c_2) = \Phi(a,c)
   -q  a_1 + p  c_2.
    \]

\begin{defn}
The restriction of $(\I_2, \Omega)$ to $Y_{\ast}$ is denoted by $(\I_{\ast}, \Omega)$
and is referred to as the
{\it Pfaffian exterior differential system of lipid bilayer membranes} satisfying the
differential relation
\[
 \Delta H = \Phi(a,c).
  \]
\end{defn}

\begin{remark}
The integral manifolds of $(\I_*, \Omega)$ are canonical prolongations of umbilic
free immersions $F : X \to \R^3$ whose mean curvature $H$ satisfies
$\Delta H = \Phi(a,c)$.

\end{remark}

\section{Involution}\label{s:PDS-involution}

\subsection{Algebraic generators}

By construction, the Pfaffian differential system $\I_*$ is generated,
as a differential ideal, by the restrictions to $Y_*$ of the 1-forms
$\alpha^j$, $\beta^a$, $\gamma^a$, $\delta^a$, $j = 1,\dots, 4; a= 1,2$.
The generators $\alpha^j, \beta^a, \gamma^a$ are expressed in terms of
the fiber coordinates  $p,q,a,c,r,a_1,c_2,p_1,q_2$ as in \eqref{alphas}
and \eqref{betas-gammas-deltas}, while the generators $\delta^1$ and
$\delta^2$ are given by
\begin{equation}\label{deltas-ell}
\begin{cases}
 \delta^1 = d a_1 - \left(\ell + r (c -a) + \Psi \right)\theta^1
           + \left( p_1 (c-a) - 2  a_1 p\right) \theta^2, \\
  \delta^2 = d c_2 + \left( q_2 (c-a) + 2 c_2 q\right) \theta^1
           + \left(\ell - r (c -a) - \Psi \right)\theta^2.\\
   \end{cases}
     \end{equation}
Observe that
\[
 (\theta^1, \theta^2, \alpha^1, \dots \alpha^4,
 \beta^1, \beta^2, \gamma^1, \gamma^2, \delta^1, \delta^2,
  d p_1, d q_2, dr, d\ell)
   \]
is a global coframe field on $Y_*$. Let
 \[
 \left( {\partial}_{\theta^1},  {\partial}_{\theta^2},
  {\partial}_{\alpha^1}, \dots,  {\partial}_{\alpha^4},
 {\partial}_{\beta^1},  {\partial}_{\beta^2},
   {\partial}_{\gamma^1},  {\partial}_{\gamma^2},
    {\partial}_{\delta^1},  {\partial}_{\delta^2},
   {\partial}_{d p_1},
    {\partial}_{d  q_2},
    {\partial}_{ dr},
    {\partial}_{d\ell}\right)
   \]
denote its dual basis on $T(Y_*)$.
According to \eqref{d-alphas-d-gammas=0},
we have
\begin{equation}\label{d-alphas-d-gammas=0*}
\begin{cases}
  d\alpha^j \equiv 0, \quad j=1,2,3,4, \\
    d\gamma^a \equiv 0, \quad a=1,2, \\
\end{cases}
\mod \{\alpha^j, \beta^a, \gamma^a, \delta^a \}.
\end{equation}
Thus, the differential ideal $\I_*$ is algebraically generated by
$\alpha^j$, $\beta^a$, $\gamma^a$, $\delta^a$, $j=1,2,3,4; a = 1,2$,
and by the exterior differential 2-forms
$d\beta^a$, $d\delta^a$, $a = 1,2$. Using \eqref{eq-Y1-1}, \eqref{betas-gammas-deltas},
\eqref{deltas-ell} and \eqref{d-alphas-d-gammas=0*}, we obtain,
modulo the algebraic ideal $\{\alpha^j, \beta^a, \gamma^a, \delta^a \}$,
\begin{equation}\label{d-betas-d-deltas} %{d-betas-d-deltas*}
\begin{cases}
  d\beta^1 \equiv  -d p_1\wedge \theta^1
        -dr \wedge \theta^2 -B^1 \theta^1 \wedge \theta^2,\\
    d\beta^2 \equiv  -dr \wedge \theta^1 -d q_2\wedge \theta^2
         -B^2 \theta^1 \wedge \theta^2,\\
 d\delta^1 \equiv  -d\ell \wedge \theta^1 - (c-a)dr\wedge\theta^1
  +(c-a)d p_1 \wedge\theta^2 -D^1\theta^1\wedge \theta^2,\\
  d\delta^2 \equiv  (c-a)d q_2 \wedge \theta^1 + d\ell\wedge\theta^2
  -(c-a)dr \wedge\theta^2 +D^2\theta^1\wedge \theta^2,\\
\end{cases}
\end{equation}
where $B^a$ and $D^a$, $a=1,2$, are real analytic functions of the fiber coordinates.

\begin{remark}
Note that if $\Psi$ is polynomial, the $B^a$ and the $D^a$ are polynomial functions
of the fiber coordinates.
This is still true of the $B^a$, but not of the $D^a$, in the event that $\Psi$ is not a
polynomial.
\end{remark}

In conclusion, we have proved that $(\I_*, \Omega)$ is algebraically generated by
$\alpha^j$, $\beta^a$, $\gamma^a$, $\delta^a$ and by

\begin{equation}\label{d-betas-d-deltas*}
\begin{cases}
  \Omega^1 =  d p_1\wedge \theta^1
        + dr \wedge \theta^2  +B^1 \theta^1 \wedge \theta^2,\\
    \Omega^2 =   dr \wedge \theta^1  +d q_2\wedge \theta^2
         + B^2 \theta^1 \wedge \theta^2,\\
 \Omega^3 =  d\ell \wedge \theta^1 +  (c-a)dr\wedge\theta^1
  - (c-a)d  p_1 \wedge\theta^2 + D^1\theta^1\wedge \theta^2,\\
  \Omega^4 =  (c-a)d q_2 \wedge \theta^1 + d\ell\wedge\theta^2
  -(c-a)dr \wedge\theta^2 +D^2\theta^1\wedge \theta^2.\\
\end{cases}
%\mod \{\alpha^j, \beta^a, \gamma^a, \delta^a \}
\end{equation}

\subsection{Polar equations}

Let $[\xi]$ be a 1-dimensional integral element,
where $\xi \in T_m Y_*$ is a tangent vector of the form
\[
 \xi = x^1  {\partial}_{\theta^1} + x^2  {\partial}_{\theta^2}
 +x^3 {\partial}_{ d p_1} + x^4  {\partial}_{ d q_2}
  +  x^5  {\partial}_{dr} + x^6  {\partial}_{d\ell},
   \quad (x^1)^2 +(x^2)^2 \neq 0.
   \]
Its {\it polar space} $H([\xi])$ is the subspace tangent to $Y_*$ defined by the
{\it polar equations}
\[
   \alpha^j = 0, \quad \beta^a= 0, \quad \gamma^a =0, \quad \delta^a =0,
   \quad i_{\xi} \Omega^j = 0, \quad j=1,\dots,4; a = 1,2.
     \]
From \eqref{d-betas-d-deltas*}, it follows that the polar equations
are linearly independent provided $c-a\neq 0$, for every 1-dimensional integral element
$[\xi]
%\in \mathcal V(\I_*, \Omega)
$.
This implies that $H([\xi])$ is the unique 2-dimensional integral element
containing $[\xi]$. Thus, by the Cartan--K\"ahler theorem we have the
following.

\begin{thm}\label{cartan-kaehler}
The Pfaffian differential system $(\I_*, \Omega)$ is in involution and its
general solutions depend on four functions in one variable.
More precisely, for every, 1-dimensional real analytic integral
manifold $\mathcal A \subset Y_*$ there exists a real analytic 2-dimensional
integral manifold $\mathcal X \subset Y_*$ passing through $\mathcal A$.
In addition, the 1-dimensional integral manifolds depend on four functions
in one variable.
\end{thm}

\begin{remark}
Note that the functional dependence of the general solutions agrees with
that of the initial data of Theorem \ref{main-thm}. In fact, besides
constants, specifying the initial data of Theorem \ref{main-thm} amounts to
the choice of a unit-speed curve $\alpha$ in $\R^3$, which depends on two
arbitrary functions in one variable,
and two functions $h$, $h^W$.
\end{remark}

\begin{remark}\label{r:4}
The 2-dimensional integral manifold $\mathcal X$ passing through $\mathcal A$
is unique, in the sense that if $\tilde {\mathcal X} \subset Y_*$ is another
real analytic 2-dimensional integral manifold containing $\mathcal A$, then
$\mathcal X$ and $\tilde {\mathcal X}$ agree in a neighborhood of $\mathcal A$.

\end{remark}

\section{The proof of Theorem \ref{main-thm}}\label{s:main-thm}

We are now ready to prove Theorem \ref{main-thm}.
We are given the Cauchy data $\alpha$, $x_0$, $W_0$, $h$, $h^W$, consisting of
a real analytic curve $\alpha : J \to \R^3$, a point $x_0\in J$,
a unit normal vector $W_0=\cos a_0 N(x_0) +\sin a_0 B(x_0)$,
and two real analytic functions
$h$, $h^W : J \to \R^3$  as in the statement of Theorem \ref{main-thm}.
Set $$s(x) =- \int_{x_0}^x{\tau(u) du + a_0},$$ so that $s : J \to \R$ is real analytic and
$s(x_0) = a_0$, and
define
\[
   m := - h - \kappa \sin s(x)
  % \left( - \int_{x_0}^x{\tau(u) du + a_0} \right)
    > 0.
  \]
Let $W : J \to \R^3$ be the unit normal vector field along $\alpha$ defined by
\[
 W = \cos s(x)
 %\left( -\int_{x_0}^x{\tau(u) du + a_0} \right)
 N
  +\sin s(x)
  %\left( -\int_{x_0}^x{\tau(u) du + a_0} \right)
  B.
    \]
Consider the positively oriented orthonormal frame field given by
\begin{equation}
   \mathcal G = (\alpha, T, W, \text{\sc J} W) : J \to \E(3),
    \end{equation}
where
\[
  \text{\sc J} W = -\sin s(x) N
    +\cos s(x) B.
     \]
Using the Frenet--Serret equations $T' = \kappa N$, $N' = -k T +\tau B$, $B' = -\tau N$
for the Frenet frame $(T,N,B)$, it is easily seen that
\[
 \frac{d\mathcal G}{dx} = \mathcal G
  \begin{pmatrix}
   0 & 0 & 0 & 0 \\
   1 & 0 & -\mathfrak p & -\mathfrak a \\
   0 & \mathfrak p & 0 & 0  \\
   0 & \mathfrak a & 0  & 0 \\
   \end{pmatrix},
     \]
where $\mathfrak p$, $\mathfrak a : J \to \R$ are the real analytic functions
\[
 \mathfrak p  = \kappa \cos s(x),
  \quad \mathfrak a = -\kappa \sin s(x).
   \]

\begin{remark}
The frame field $(T, W, \text{\sc J} W)$ is a relatively parallel adapted frame along
the curve $\alpha$ in the sense of Bishop \cite{Bishop}.
\end{remark}

Next, if we set
\begin{equation}\label{gothic-p-q-r-l}
\begin{cases}
\mathfrak c = \mathfrak a - 2m ,\\
    \mathfrak q = -\frac{1}{\mathfrak c - \mathfrak a}\frac{d\mathfrak c}{dx},\\
    {\mathfrak a}_1 = \frac{d\mathfrak a}{dx}, \\
      {\mathfrak c}_2 = 2h^{W}
         + \mathfrak p (\mathfrak c - \mathfrak a), \\
     {\mathfrak p}_1 = \frac{d\mathfrak p}{dx}, \\
       {\mathfrak q}_2  = -\frac{1}{\mathfrak c - \mathfrak a}
         \left(\frac{d {\mathfrak c}_2}{dx}
              + 2 {\mathfrak c}_2\mathfrak q\right), \\
        \mathfrak r = \frac{d{\mathfrak q}}{dx}
             + \frac{1}{2}(\mathfrak a \mathfrak c + \mathfrak p^2
          + \mathfrak q^2), \\
     \mathfrak l = \frac{d^2{\mathfrak a}}{dx^2}
              - \mathfrak r(\mathfrak c - \mathfrak a) - \Psi(\mathfrak p,
      \mathfrak q, {\mathfrak a}_1, {\mathfrak c}_2),
\end{cases}
\end{equation}
then
\begin{equation}\label{def-mathcal-A}
  \mathcal A : J \to Y_*, \, x \mapsto \left(\mathcal G, \mathfrak p, \mathfrak q,
   \mathfrak a, \mathfrak c, {\mathfrak p}_1,  {\mathfrak q}_2,
    \mathfrak r,
     {\mathfrak a}_1,  {\mathfrak c}_2, \mathfrak l\right)\vert_x
     \end{equation}
is a 1-dimensional integral manifold of $\I_*$ such that
 \[
   \theta^1 = dx, \quad \theta^2 = 0,
    \]
    defined by
\begin{equation}\label{gothic-p-q-r-l-A}
\begin{cases}
  p \circ \mathcal A = \mathfrak p ,\quad
  q \circ \mathcal A = \mathfrak q ,\\
    a \circ \mathcal A = \mathfrak a,\quad
    c \circ \mathcal A = \mathfrak c,\\
   p_1 \circ \mathcal A =  {\mathfrak p}_1, \quad
   q_2 \circ \mathcal A =  {\mathfrak q}_2, \\
   r \circ \mathcal A = \mathfrak r, \\
   a_1 \circ \mathcal A =   {\mathfrak a}_1, \quad
   c_2 \circ \mathcal A =   {\mathfrak c}_2, \\
      l\circ \mathcal A = \mathfrak l.
\end{cases}
%\quad (a = 1,2)
\end{equation}

\begin{defn}
We call $\mathfrak U : = \text{Im}\, \mathcal A\subset Y_*$ the {\it canonical
1-dimensional integral manifold} defined by the {\it Cauchy data}
$\left( \alpha, x_0, W_0, h, h^{W}\right)$.
\end{defn}

For a set of Cauchy data $\left( \alpha, x_0, W_0, h, h^{W}\right)$,
Theorem \ref{cartan-kaehler} and Remark \ref{r:4} imply that there exists
a unique real analytic integral manifold $\mathcal X \subset Y_*$ of
$(\I_*, \Omega)$ such that $\mathfrak U \subset \mathcal X$, where
$\mathfrak U$ is the integral curve defined by
$\left( \alpha, x_0, W_0, h, h^{W}\right)$.
On $\mathcal X$ we consider the orientation defined by the 2-form
$\Omega = \theta^1\wedge \theta^2$. The map
\[
 F : \mathcal X \ni (P, A, p, q, a, c, p_1, q_2, r, a_1, c_2, \ell) \mapsto P\in \R^3
  \]
is a real analytic immersion and, by construction, its prolongation $F_{(2)}$
coincides with the inclusion map $\iota : \mathcal X \to Y_*$. According
to Remarks 1 and 2, we infer that
\begin{itemize}
\item $\mathcal X \ni (P, A, p, q, a, c, p_1, q_2, r, a_1, c_2, \ell) \mapsto (P,A)\in \E(3)$
is a principal frame field along $F$;

\item $F$ satisfies $\Delta H = \Phi(a,c)$;

\item $(\theta^1, \theta^2)$ is a positively oriented orthonormal principal coframe
on $\mathcal X$ along the immersion $F$;

\item $a, c : \mathcal X \to \R$ are the principal curvatures of $F$ and
\[
 dH = \frac{1}{2}(da + dc) = \frac{1}{2}(a_1 -q(c-a))\theta^1
     + \frac{1}{2}( c_2 -p(c-a))\theta^2.
    \]

\end{itemize}

Since $\mathcal X$ satisfies the initial condition $\mathfrak U \subset \mathcal X$
and $\mathcal A^*(\theta^2) = 0$, we can state the following.

\begin{lemma}\label{l:4}
$\mathcal A : J \to \mathcal X$ is a curvature line of $F$ such that
\[
   F \circ \mathcal A = \alpha, \quad A\circ \mathcal A = \mathcal G.
   \]
\end{lemma}

In particular, we have
\[
 F_*(T_{\mathcal A(x_0)}(\mathcal X) = \text{span}\,
(A_1(\alpha(x_0)), A_2(\alpha(x_0)) = \text{span}\,(T(x_0), W(x_0)).
\]
From \eqref{betas-gammas-deltas}, \eqref{def-mathcal-A} and \eqref{gothic-p-q-r-l-A},
we obtain
\begin{equation}\label{H-comp-mathcal-A}
 H\circ \mathcal A = \frac{1}{2}(\mathfrak a + \mathfrak c) = h
  \end{equation}
and
\begin{equation}\label{dH-mathcal-A}
 dH\vert_{\mathcal A(x_0)} \equiv
 \frac{1}{2}( {\mathfrak c}_2
   -\mathfrak p(\mathfrak c-\mathfrak a))\theta^2\vert_{\mathcal A(x_0)}
=  h^{W} \theta^2\vert_{\mathcal A(x_0)}, \mod \theta^1\vert_{\mathcal A(x_0)}.
  \end{equation}

Let $X_1 ,X_2$ denote the frame field
dual to the coframe field $\theta^1,\theta^2$ on the integral manifold $\mathcal X$
and let $\Theta$ be the local flow generated by the nowhere vanishing vector field
$X_2$. Then $\Theta$ is a real analytic map $\Theta : \mathcal U \to \mathcal X$
defined on an open neighborhood $\mathcal U\subset \mathcal X \times \R$ and the set
\[
 \Sigma := \left\{(x, y) \in J\times \R \,\, | \,\, (\mathcal A (x), y) \in \mathcal U \right\}
  \]
is an open neighborhood of $J \times \{0\}\subset \R^2$.
Using the immersion $F$ and the flow $\Theta$, define the map
$f : \Sigma \to \R^3$ by
\[
   f(x,y) = F\left(\Theta(\mathcal A(x),y)\right).
  \]
The map $f: \Sigma \to \R^3$ is a real analytic immersion such that
\[
  f(x,0) = F\left(\Theta(\mathcal A(x), 0)\right) = \alpha (x), \, \forall \, x\in J.
   \]
By construction, $f$ is a re-parametrization of $F$ and hence satisfies the
differential relation \eqref{shape-eq}. Lemma \ref{l:4} and the equations
\eqref{H-comp-mathcal-A} and \eqref{dH-mathcal-A} imply that $f$ satisfies
the required conditions (3), (4) and (5) of Theorem \ref{main-thm}.
If $\hat f : \hat \Sigma \to \R^3$ is any other immersion satisfying the same
conditions, then its
prolongation $\hat f_{(2)} : \hat \Sigma \to Y_*$ is an integral manifold of
$(\I_2,\Omega)$ passing through $\mathfrak U$. Then, by the uniqueness part of
Cartan--K\"ahler theorem, it follows that
$f_{(2)}(\Sigma\cap\hat \Sigma) =\hat f_{(2)}(\Sigma\cap\hat \Sigma)$.
This concludes the proof of Theorem \ref{main-thm}.

\section{Examples}\label{s:ex}

In this section we illustrate Theorem \ref{main-thm} in the case of cylindrical
equilibrium configurations \cite{OY1996, VDM}. We will show how such solution
surfaces are related to congruence curves, that is, plane curves that move
without changing their shape when their curvature evolves according to the
modified KdV equation \cite{GP, NSW, M-cm}.
\vskip0.1cm

In $\R^3$, with Cartesian coordinates $x_1, x_2, x_3$,
let $\Gamma\subset \R^3$ be an embedded simple closed curve contained in the
coordinate $x_1x_2$-plane, oriented by the third vector of the canonical basis, $e_3$.
%, where $(e_1,e_2,e_3)$ is the canonical basis.
Consider a unit-speed parametrization $\alpha : \R\to \R^2 \subset\R^3$ of $\Gamma$, which in turn
determines an orientation along $\Gamma$. We let $\epsilon=\pm 1$ according to whether the
orientation induced by $\alpha$ is counterclockwise or clockwise.
Let $J\alpha'$ be the
unit normal vector field along $\alpha$ obtained by counterclockwise rotating the unit
tangent vector field $\alpha'$ by an angle $\pi/2$, and
let $\kappa = \alpha''\cdot J\alpha' :\R \to \R$ be the signed curvature of $\alpha$.
Let $S \subset \R^3$ be the cylinder
with directrix curve $\Gamma$ and generating lines parallel to $e_3$.
On $S$ we put the orientation determined by the
outward unit normal vector field.
Then, $f(x,y)=\alpha(x)+\epsilon ye_3$ is a parametric equation of $S$ compatible with the
given orientation. The principal curvatures are $a=-\epsilon \kappa$ and $c=0$, so that the
Gaussian curvature vanishes identically and $H=-\epsilon \kappa/2$.
It then follows that $S$ satisfies the equation
$\Delta H=\Phi(a,c)$ if and only if the signed curvature $\kappa$ is a solution of the second
order ordinary differential equation
\begin{equation}\label{eq1}
    \kappa''=-2\epsilon\Phi(-\epsilon\kappa,0).
    \end{equation}
Note that the cylinder $S$ has no umbilic points if and only if its directrix curve $\Gamma$
is strictly convex (i.e., $\kappa$ is either strictly positive or strictly negative according to
whether $\Gamma$ is oriented counterclockwise or clockwise). In this case, the answer to the
geometric Cauchy problem provided by Theorem \ref{main-thm} is the following.
Suppose we are given
\begin{itemize}
\item a convex simple closed plane curve $\Gamma\subset \R^2(x_1,x_2)\subset \R^3$, with signed curvature $\kappa$ satisfying \eqref{eq1};
\item a point $\alpha(x_0)\in \Gamma$ and the unit normal vector $W_0 =-\epsilon e_3$ (it corresponds to the value $-\pi/2$ of the constant $a_0$);
\item $h=-\epsilon \kappa/2$ and $h^W=0$,
\end{itemize}
as Cauchy data. Then, the integral inequality \eqref{kt-ineq} is fulfilled and the cylinder
$S$ is the unique surface satisfying $\Delta H  = \Phi(a,c)$ determined by the given initial data.

\vskip0.1cm

\subsection{Cylindrical membranes} We now focus on the membrane shape equation
\eqref{shape-eq-alt1}. In this case, \eqref{eq1} takes the form
\begin{equation}\label{eq2}
  \kappa''+\frac{1}{2}\kappa^3 - \frac{2\lambda + kc_0^2}{2k}\,\kappa - \frac{\epsilon p}{k}=0.
   \end{equation}
Putting $v=(2\lambda+kc_0^2)/2k$ and differentiating (\ref{eq2}), we obtain
\[
  \kappa'''+ \frac{3}{2}\kappa^2\kappa'-v \kappa'=0,
   \]
which implies that the function $(x,t)\in \R^2 \mapsto \kappa(x+vt)\in \R$ is a traveling wave solution
of the modified KdV (mKdV) equation
\cite{GP}. From a geometrical point of view this is equivalent to saying
that $\Gamma$ is a Goldstein--Petrich contour, that is, a simple
closed congruence curve of the mKdV flow (second Goldstein--Petrich flow) \cite{M-cm}.
Equation (\ref{eq2}) admits a first integral, namely,
multiplication of \eqref{eq2} by $\kappa'$ and integration yields
\begin{equation}\label{eq3}
 (\kappa')^2+\frac{1}{4}(\kappa^4+w_2\kappa^2+w_1\kappa+w_0)=0,
  \end{equation}
where $w_0\in \R$ is a constant of integration, $w_2=-2(2\lambda +kc_0^2)/k$ and $w_1=-8\epsilon p/k$.
If the pressure $p$ vanishes identically, then $S$ is a Willmore cylinder and $\Gamma$
is a planar elastica.
It is well known that all closed elastic planar curves are lemniscates and
therefore all of them possess a point of self-intersection (see \cite{Br-Gr}, for instance).
Discarding this case and by possibly rescaling $\Gamma$ by the similarity factor $r=(w_1)^{1/3}$,
we may assume $w_1=1$. The solutions of \eqref{eq3} can be explicitly computed in terms of
Jacobi's elliptic functions. The closedness conditions and the embeddedness of plane curves
whose signed curvature satisfies \eqref{eq3} have been investigated independently, and
in different contexts, by Vassilev--Djondjorov--Mladenov \cite{VDM} and Musso \cite{M-cm}.
We briefly summarize the results obtained in \cite{M-cm}.
Suppose that the polynomial $t^4+w_2t^2+t+w_0$ has two distinct real roots and two
complex conjugate roots. Then, the constants $w_0$ and $w_2$ can be written as

\begin{equation}\label{c}
  w_{0}=\frac{(1+4\varsigma^3\varrho^2)(1+4\varsigma^3(\varrho^2-1))}{16\varsigma^4},\quad
    w_{2}=-\frac{1}{2\varsigma^2}+\varsigma(2\varrho^2-1),
    \end{equation}
where $\varsigma$ and $\varrho$ are two parameters such that $\varsigma<0$ and $\varrho\in(-1,1)$. Let
\begin{equation}\label{gm}
 g=-\frac{1}{2\varsigma}\left(1+\varsigma^6+\varsigma^3(4\varrho^2-2) \right)^{1/4},\quad
  m=\frac{1}{2}+\frac{\varsigma^3(1-2\varrho^2)-1}{2\left(1+\varsigma^6+\varsigma^3(4\varrho^2-2)\right)^{1/2}},
   \end{equation}
and define $\alpha_1$, $\alpha_2$, $\beta_1$ and $\beta_2$ by posing
\begin{equation}\label{alpha}
\begin{cases}
 \alpha_{1}=A_{1}-A_{2},\quad &\alpha_{2}=-(A_{1}+A_{2}),\\
  \beta_{1}=B_{1}-B_{2},\quad &\beta_{2}=-(B_{1}+B_{2}),
   \end{cases}
\end{equation}
where $A_1,A_2,B_1$ and $B_2$ are given by
\[
 \begin{cases}
  A_{1}=\frac{1}{2\varsigma^2}\sqrt{1-\varsigma^3+2\varrho(-\varsigma)^{3/2}}(1-2\varrho(-\varsigma)^{3/2}),\\
  A_{2}= \frac{1}{2\varsigma^2}\sqrt{1-\varsigma^3-2\varrho(-\varsigma)^{3/2}}(1+2\varrho(-\varsigma)^{3/2}),\\
  B_{1}=\frac{1}{\varsigma}\sqrt{1-\varsigma^3+2\varrho(-\varsigma)^{3/2}},\\
  B_{2}=\frac{1}{p\varsigma}\sqrt{1-\varsigma^3-2\varrho(-\varsigma)^{3/2}}.
   \end{cases}
   \]
It is now a computational matter to check that
\begin{equation}\label{curvature}
  \kappa_{\varsigma,\varrho}(s) = \frac{\alpha_{1}\mathrm{cn}(g s|m)+\alpha_{2}}
    {\beta_{1}\mathrm{cn}(gs|m)+\beta_{2}}
     \end{equation}
is a periodic solution of \eqref{eq3}. The period of $\kappa_{\varsigma,\varrho}$ is the complete elliptic integral
\begin{equation}\label{period}
 \omega_{\varsigma,\varrho}=\frac{4}{g}\int_0^{\pi/2}\frac{d\vartheta}{\sqrt{1-m\sin^2 \vartheta}}.
  \end{equation}
We consider the angular function
\[
  \theta_{\varsigma,\varrho}(s):=\int_0^s \kappa_{\varsigma,\varrho}(u)du.
  \]
Such a function can possibly be explicitly evaluated in terms of incomplete
elliptic integrals of the third kind. Setting
\[
  \alpha_{\varsigma,\varrho}:=2\left((2\kappa_{\varsigma,\varrho}+w_2)\cos \theta_{\varsigma,\varrho}- 4\kappa_{\varsigma,\varrho}' \sin\theta_{\varsigma,\varrho},
     (2\kappa_{\varsigma,\varrho}+w_2)\sin\theta_{\varsigma,\varrho}
     -4\kappa_{\varsigma,\varrho}'\cos\theta_{\varsigma,\varrho},0\right),
    \]
we obtain a unit-speed clockwise parametrization of an immersed curve $\Gamma_{\varsigma,\varrho}$, with signed curvature $\kappa_{\varsigma,\varrho}$. In addition, $\alpha_{\varsigma,\varrho}$ is periodic if and only if
\[
 \Lambda_{\varsigma,\varrho}=\frac{1}{2\pi}\int_0^{\omega_{\varsigma,\varrho}} \kappa_{\varsigma,\varrho}(u)du =
  \frac{\mu}{\upsilon}\in \mathbb{Q},
   \]
where $\mu,\upsilon\in \Z$ are relatively prime integers, with $\mu \ge 0$.

\begin{remark}
If $\varrho=0$, then $\kappa_{\varsigma,0}=1/2\varsigma$ and the corresponding curve is a circle of radius $2|\varsigma|$.
If $\varrho \neq 0$, the integers $\mu$ and $\upsilon$ have the following geometrical meaning:
$\mu$ is the turning number and $|\upsilon|$ is the order of the symmetry group of the
immersed curve. In particular, for a simple curve, $\mu=1$.
\end{remark}

In \cite{M-cm}, it is proved that, for every positive integer $\upsilon >1$, there exists
$\varrho_{\upsilon}\in (0,1)$ and a real-analytic map $\phi_{\upsilon}:[0,\varrho_{\upsilon})\to \R$,
such that $\alpha_{\phi_{\upsilon}(\varrho),\varrho}$ is strictly convex. When $\varrho\in (0,\varrho_{\upsilon})$,
the corresponding curve has a non-trivial symmetry group generated by a rotation of an angle
$2\pi/\upsilon$ around the $x_3$-axis. If $\varrho=0$, then $\phi_{\upsilon}(\varrho)=-(\upsilon^2-1)^{1/3}$,
and we get a circle of radius $2(\upsilon^2-1)^{1/3}$.
This shows that, for every positive integer $\upsilon>1$, there exists a 1-parameter
family $S_{\upsilon,\varrho}$, $\varrho\in [0,\varrho_{\upsilon})$, of embedded cylindrical membranes
which are invariant under the subgroup generated by the rotation of an angle $2\pi/\upsilon$
around the $x_3$-axis. If $\varrho=0$, then $S_{\upsilon,0}$ is a round cylinder.
We would like to stress that $\varrho_{\upsilon}$ and the function $\phi_{\upsilon}$ can be
evaluated by numerical methods, so that all such cylindrical solution surfaces can be
effectively determined. The codes for such computations can be found in \cite{M-cm}.

\begin{figure}[ht]
\begin{center}
\includegraphics[height=5cm,width=5cm]{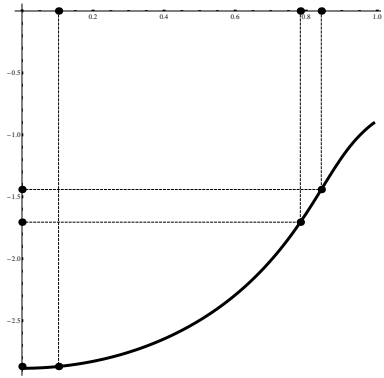}
\caption{The graph of the function $\phi_5$ and the ``separating values'' $\varrho_5$, $\varrho_{1,5}$ and $\varrho_{2,5}$.}\label{FIG1}
\end{center}
\end{figure}

\begin{remark}
Numerical experiments exhibit some additional geometric
properties of these families of cylindrical membranes. First, the function $\phi_{\upsilon}$
can be defined on the whole interval $[0,1)$ and, if $\lfloor \upsilon/2\rfloor$ denotes
the integer part of $\upsilon/2$, there exists a strictly increasing sequence
$\{\varrho_{j,\upsilon}\}_{j=1,\dots,\lfloor \upsilon /2\rfloor} \subset (\varrho_{\upsilon},1)$, such that:
 \begin{itemize}
\item  the curve $\Gamma_{\phi_{\upsilon}(\varrho_{\upsilon}),\varrho_{\upsilon}}$ is convex, but not strictly convex,
    and has exactly $\upsilon$ inflexion points corresponding to
    $\alpha_{\phi_{\upsilon}(\varrho),\varrho}(j\omega_{\phi_{\upsilon}(\varrho),\varrho}/2)$, $j=1,\dots,\upsilon$;
\item if $\varrho>\varrho_{\upsilon}$,
    $\Gamma_{\phi_{\upsilon}(\varrho),\varrho}$ is not convex
    and has $2\upsilon$ inflection points;
\item if $\varrho\in [\varrho_{\upsilon},\varrho_{1,\upsilon})$, $\Gamma_{\phi_{\upsilon}(\varrho),\varrho}$ is a simple star-like curve;
\item for every $j=1,\dots,\lfloor\upsilon/2\rfloor$,   $\Gamma_{\phi_{\upsilon}(\varrho_{j,\upsilon}),\varrho_{j,\upsilon}}$ has $\upsilon(2j-1-\delta_{2j,n})$ points of self-intersections,
    of which $(1 - \delta_{2j,n}/2)\upsilon$ are non-transversal;
\item if $\varrho\in (\varrho_{j,\upsilon},\varrho_{j+1,\upsilon})$, $\Gamma_{\phi_{\upsilon}(\varrho),\varrho}$ has $2j\upsilon$ transversal self-intersections;
\item if $\varrho\in (\varrho_{\lfloor\upsilon/2\rfloor,\upsilon},1)$, $\Gamma_{\phi_{\upsilon}(\varrho),\varrho}$ has $\upsilon(\upsilon-1)$ transversal self-intersections.
\end{itemize}
\end{remark}

\begin{figure}[ht]
\begin{center}
\includegraphics[height=4cm,width=8cm]{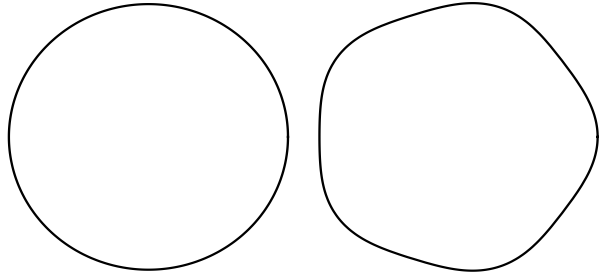}
\caption{The curves $\Gamma_{\phi_5(0),0}$ and $\Gamma_{\phi_5(\varrho),\varrho}$, $\varrho=0.08$.}\label{FIG2}
\end{center}
\end{figure}

We would like to point out that the explicit determination of the membrane surfaces
which are invariant under a 1-parameter subgroup of rigid motions is,
in its full generality, an open problem.
Rotationally invariant solutions have been considered by Capovilla--Guven--Rojas \cite{CGR},
but little is known about the more general class of helicoidal solution surfaces \cite{Tek2007}.
Hopefully, our approach may be useful to tackle this problem; in fact,
with an additional invariance hypothesis,
our exterior differential system is not
in involution anymore, but it is to expect that a suitable prolongation of it
would be completely integrable, in the sense of Frobenius.
In principle, the problem could then be reduced to a system of ordinary
differential equations.

\begin{figure}[ht]
\begin{center}
\includegraphics[height=4cm,width=8cm]{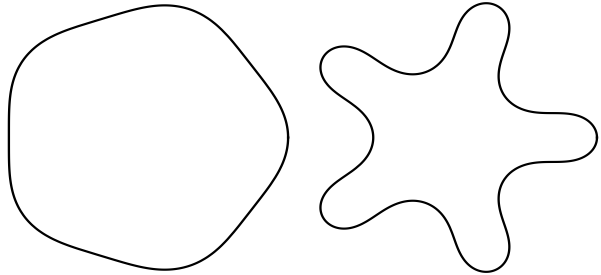}
\caption{The curves $\Gamma_{\phi_5(\varrho_5),\varrho_5}$ and $\Gamma_{\phi_5(\varrho),\varrho}$, $\varrho=0.6\in (\varrho_5,\varrho_{1,5})$.}\label{FIG3}
\end{center}
\end{figure}

\subsection{Numerical computations and visualization} We now illustrate the phenomenology
of the 1-parameter family of cylindrical membranes with a five-fold symmetry.
The numerical computations and the visualization are based on the codes given in \cite{M-cm}.
\vskip0.1cm

$\bullet$ Figure \ref{FIG1} shows the graph of the $\phi_5$-function and
the ``separating values'' $\varrho_5$, $\varrho_{1,5}$ and $\varrho_{2,5}$ of the parameter $\varrho\in [0,1)$.
The approximate values are: $\varrho_5 \approx 0.103$, $\varrho_{1,5} \approx 0.783468$ and
$\varrho_{2,5}\approx 0.84245$.
\vskip0.1cm

\begin{figure}[ht]
\begin{center}
\includegraphics[height=4cm,width=8cm]{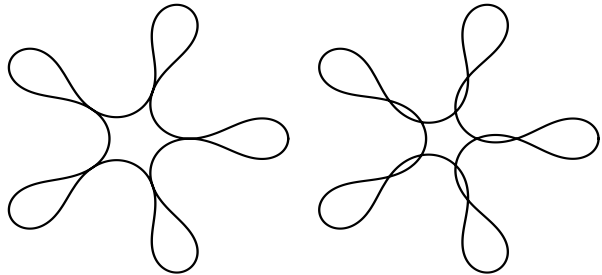}
\caption{The curves $\Gamma_{\phi_5(\varrho_{1,5}),\varrho_{1,5}}$ and $\Gamma_{\phi_5(\varrho),\varrho}$, $\varrho=0.8\in (\varrho_{1,5},\varrho_{2,5})$.}\label{FIG4}
\end{center}
\end{figure}

$\bullet$ When $\varrho=0$, the directrix is a circle of radius $4\sqrt[3]{3}$.
If the parameter $\varrho\in(0,\varrho_5)$, the directrix curve is a strictly convex
star-like contour (see Figure \ref{FIG2}).
\vskip0.1cm

$\bullet$ If $\varrho=\varrho_5$, the curve is convex, but not strictly convex. Its
inflection points are located at $\alpha_{\phi_5(\varrho),\varrho}(r\omega_{\phi_5(\varrho),\varrho}/2)$,
$r=1,\dots,5$. When $\varrho\in (\varrho_5,\varrho_{1,5})$, the directrix is not convex, with
ten inflexion points (see Figure \ref{FIG3}).
\vskip0.1cm

$\bullet$ If $\varrho=\varrho_{1,5}$, the directrix has five non-transversal self-intersections.
When $\varrho\in (\varrho_{1,5}, \varrho_{2,5})$,
the curve has ten self-intersections (see Figure \ref{FIG4}).
\vskip0.1cm

\begin{figure}[ht]
\begin{center}
\includegraphics[height=4cm,width=8cm]{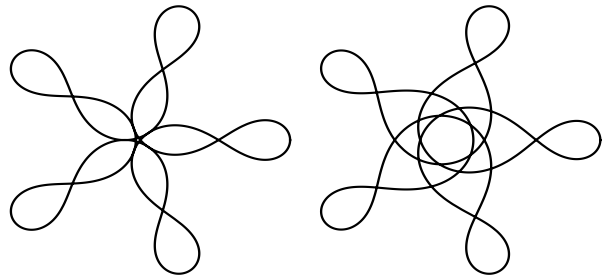}
\caption{The curves $\Gamma_{\phi_5(\varrho_{2,5}),\varrho_{2,5}}$ and $\Gamma_{\phi_5(\varrho),\varrho}$,
$\varrho=0.9\in (\varrho_{2,5},1)$.}\label{FIG5}
\end{center}
\end{figure}

$\bullet$ Finally, if $\varrho=\varrho_{2,5}$, the curve has fifteen self-intersections,
five of which are non-transversal. If $\varrho\in (\varrho_{2,5},1)$, the curve
has twenty transversal self-intersections (see Figure \ref{FIG5}).

\begin{figure}[ht]
\begin{center}
\includegraphics[height=5cm,width=5cm]{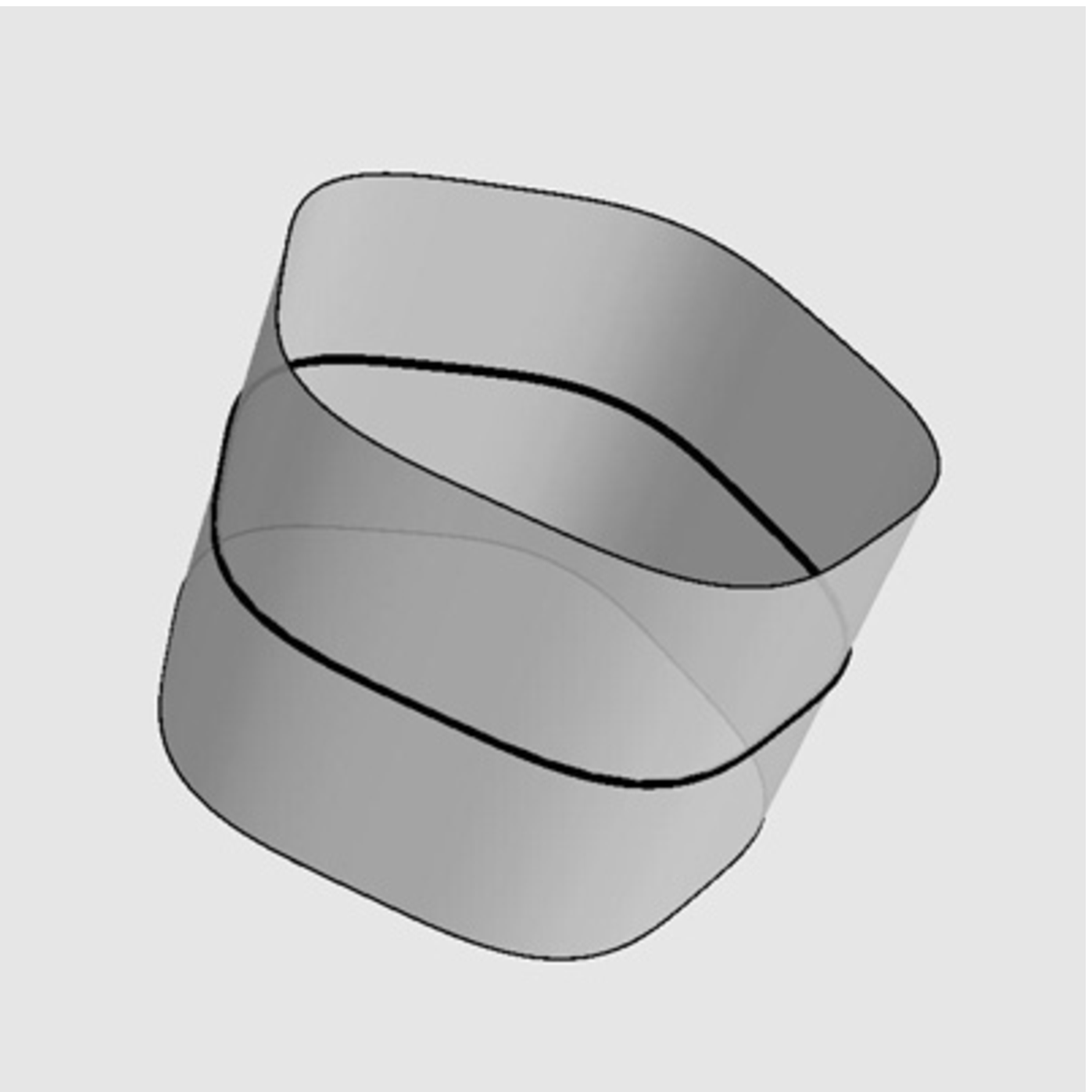}
\includegraphics[height=5cm,width=5cm]{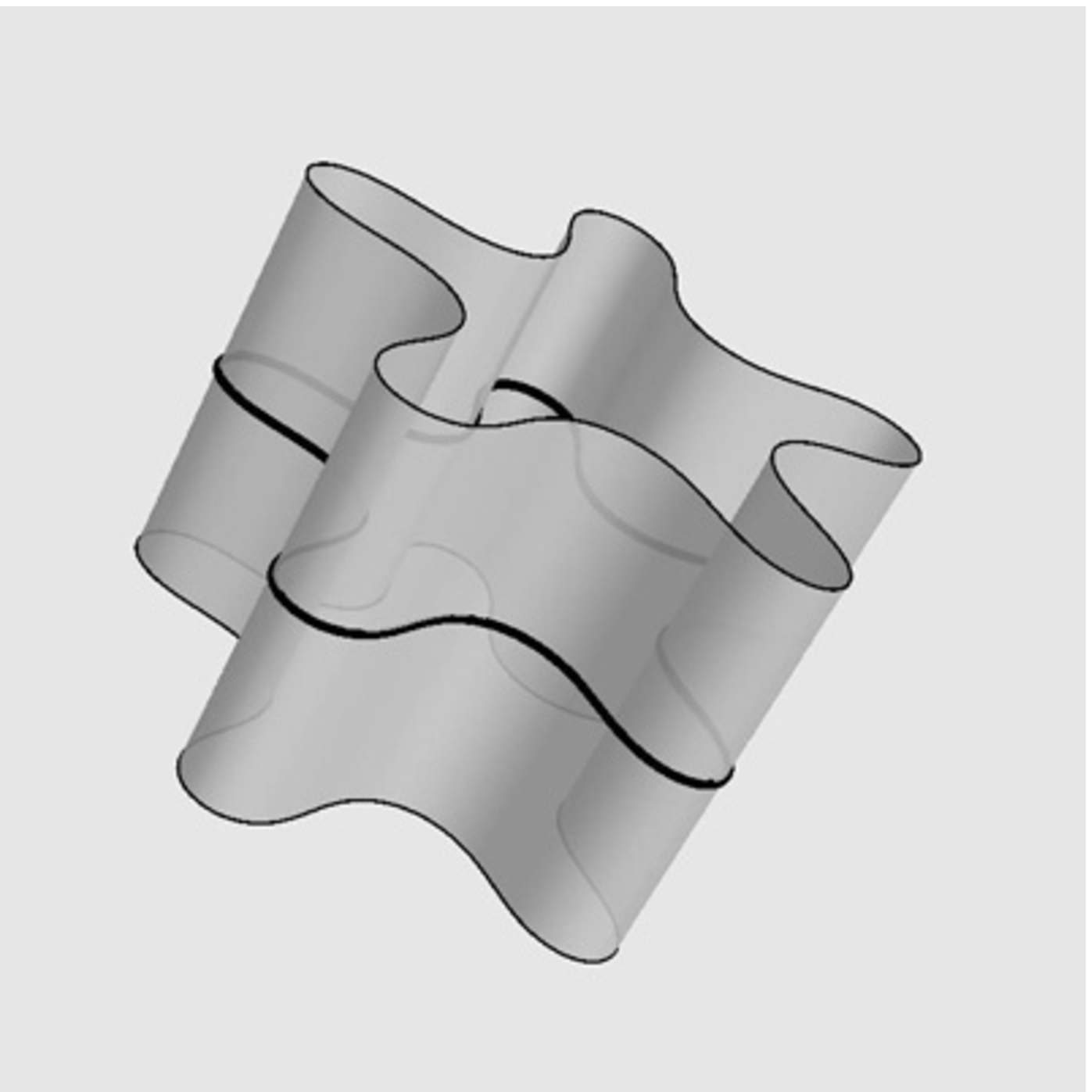}
\caption{Two cylindrical membranes of the family with a five-fold symmetry corresponding
to the values $\varrho=0.10300000000000001$ and $\varrho=0.6$ of the parameter $\varrho$.}\label{FIG6}
\end{center}
\end{figure}

\begin{figure}[ht]
\begin{center}
\includegraphics[height=5cm,width=5cm]{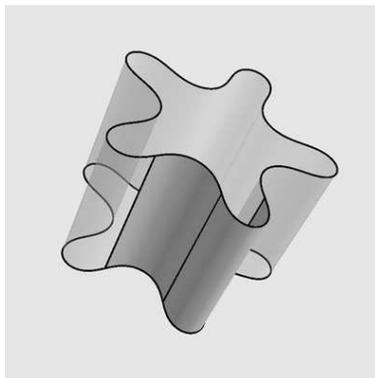}
\caption{Two strips of the cylinder with $\upsilon=5$ and $\varrho=0.6$ which can be deduced from Theorem \ref{main-thm}.}\label{FIG8}
\end{center}
\end{figure}

\noindent Our result allow to deduce the existence of global solutions only when $\varrho\in [0,\varrho_5)$
(i.e., if the directrix is strictly convex).
Figure \ref{FIG6} reproduces the cylindrical membranes of the family with a five-fold symmetry
corresponding to the values $\varrho=0.10300000000000001\in [0,\varrho_5)$ and $\varrho=0.6\in (\varrho_5,\varrho_{1,5})$.
The global existence of the cylinder on the
left
can be inferred from Theorem \ref{main-thm}.
When the directrix has inflection points, then using Theorem \ref{main-thm}, we can only foresee
the vertical strips of the cylinder between the generating lines passing through two
consecutive inflection points of the directrix curve (see Figure \ref{FIG8}).

\bibliographystyle{amsalpha}

\end{document}